\documentclass[reqno]{amsart}
\usepackage{amscd}
\usepackage{amssymb}
\usepackage{amsmath}
\usepackage{amsthm}
\usepackage{latexsym}
\usepackage{upref}
\usepackage{epsfig}
\usepackage{mathrsfs}
\usepackage{float}
\usepackage[colorlinks,urlcolor=black]{hyperref}
\usepackage{graphics,graphicx,color}
\usepackage{floatflt}
\usepackage{indentfirst}
\usepackage[cmtip,arrow]{xy}
\usepackage{pb-diagram,pb-xy}
\usepackage{cite}

\newtheorem{thm}{Theorem}[section]
\newtheorem{prop}[thm]{Proposition}
\newtheorem{cor}[thm]{Corollary}
\newtheorem{lem}[thm]{Lemma}
\theoremstyle{definition}

\newtheorem{exa}[thm]{Example}
\numberwithin{equation}{section}

\begin{document}
\title{A Characterization of trace zero bisymmetric nonnegative $5\times 5$ matrices}

\author{Somchai Somphotphisut}
\address{Department of Mathematics and Computer Science, Faculty of Science,
Chulalongkorn University, Phyathai Road, Patumwan, Bangkok 10330,
Thailand} \email{somchai.so@student.chula.ac.th \textrm{and
}kwiboonton@gmail.com}

\author{Keng Wiboonton}
\address{}
\email{}

\begin{abstract}
Let $\lambda_1 \geq \lambda_2 \geq \lambda_3 \geq \lambda_4 \geq \lambda_5 \geq -\lambda_1$ be real numbers such that $\sum_{i=1}^5 \lambda_i =0$. In \cite{oren}, O. Spector prove that a necessary and sufficient condition for $\lambda_1, \lambda_2, \lambda_3, \lambda_4, \lambda_5$ to be the eigenvalues of a traceless symmetric nonnegative $5 \times 5$ matrix is ``$\lambda_2+\lambda_5<0$ and $\sum_{i=1}^5 \lambda_{i}^{3} \geq 0"$. In this article, we show that this condition is also a necessary and sufficient condition for $\lambda_1, \lambda_2, \lambda_3, \lambda_4, \lambda_5$ to be the spectrum of a traceless bisymmetric nonnegative $5 \times 5$ matrix.
\end{abstract}
\maketitle
\section{Introduction}
\medskip

The problem of finding necessary and sufficient conditions for $n$ complex numbers to be the spectra of a nonnegative matrix is known as the \textit{nonnegative inverse eigenvalues problem} (NIEP). This problem was first studied by Suleimanova in \cite{sul} and then several authors had been extensively studied this problem. Several results for this problem can be found in \cite{borobia, Kel, rado, salz, soto1, soto2, sul, loewy2, wuwen}. If we consider this problem only for the list of $n$ real numbers, then the problem will be called the \textit{real nonnegative inverse eigenvalue problem} (RNIEP).

The problem of finding the necessary and sufficient conditions on the list of $n$ real numbers to be a spectrum of a symmetric (bisymmetric) nonnegative matrix is called the \textit{symmetric} (\textit{bisymmetric}) \textit{nonnegative inverse eigenvalue problem} (SNIEP (BNIEP)). In \cite{julio}, A. I. Julio and R. L. Soto give sufficient conditions for the BNIEP. This article discusses the BNIEP for $5 \times 5$ traceless matrices.

Recall that a square matrix $A$ is \textit{symmetric}, \textit{persymmetric} or \textit{centrosymmetric} if $A = A^T$, $AJ = JA^T$ or $AJ = JA$, respectively, where $J$ is the reverse identity matrix, i.e., $J=\begin{bmatrix}
 0 &&1\\
 &\reflectbox{$\ddots$}&\\
1&& 0
\end{bmatrix}$ (note that $J^2 = J$). A square matrix possessing any two of these conditions possesses the third, and we call it \textit{bisymmetric}. Many results on bisymmetric matrices were discussed by A. Cantoni and P. Butler in \cite{can}.






\medskip
The NIEP, RNIEP, SNIEP and BNIEP have been unsolved for general $n$. The NIEP for $n \leq 3$ and the RNIEP for $n \leq 4$ have been solved by Loewy and London \cite{loewy2}. In \cite{wuwen}, Wuwen showed that the RNIEP and the SNIEP are equivalent when $n \leq 4$. The fact that the RNIEP and the SNIEP are different for $n > 4$ was proved by Johnson et al. in \cite{johnson}. The NIEP for $4 \times 4$ traceless matrices was solved by Reams in \cite{reams} and $5\times5$ traceless matrices was solved by Laffey and Meehan in \cite{laffey2}. In \cite{oren}, O. Spector showed that a necessary and sufficient condition for a list of $5$ real numbers $\lambda_1 \geq \lambda_2 \geq \lambda_3 \geq \lambda_4 \geq \lambda_5$ such that $\lambda_5 \geq -\lambda_1$ and $\sum_{i=1}^5 \lambda_i =0$ to be the eigenvalues of a $5 \times 5$ traceless symmetric nonnegative matrix is $\lambda_2 +\lambda_5 \leq 0$ and $\sum_{i=1}^5 \lambda_{i}^3 \geq 0$. In this article, we prove that the above condition is also a necessary and sufficient condition of the BNIEP for traceless $5 \times 5$ matrices. We organize our paper as follows. In Section $2$, we give sufficient conditions of the BNIEP for $5 \times 5$ matrices. This result is Proposition \ref{prop5by5}. Then in Section $3$, we characterize all bisymmetric nonnegative traceless $5 \times 5$ matrices. Our main result of this characterization is Theorem \ref{5by5unknowcase1}.
 
\section{Sufficient conditions of the BNIEP for $5 \times 5$ matrices}
\medskip
Let $\lambda_1 \geq \lambda_2 \geq \lambda_3 \geq \lambda_4 \geq \lambda_5$ be real numbers. If $\lambda_1, \lambda_2, \lambda_3, \lambda_4, \lambda_5$ is the eigenvalues of a $5 \times 5$ nonnegative matrix, then $\sum_{i=1}^5 \lambda_i \geq 0$ and by the Perron-Frobenius Theorem, $\lambda_1 \geq \vert \lambda_5 \vert$. Since we consider the BNIEP for $5 \times 5$ matrices, we from now on assume that $\lambda_1 \geq \lambda_2 \geq \lambda_3 \geq \lambda_4 \geq \lambda_5 \geq -\lambda_1$ and $\sum_{i=1}^5 \lambda_i \geq 0$.

\medskip
If $\lambda_1 \geq \lambda_2 \geq \lambda_3 \geq \lambda_4 \geq \lambda_5 \geq 0$ or $\lambda_1 \geq \lambda_2 \geq \lambda_3 \geq \lambda_4 \geq 0 > \lambda_5$ then a $5 \times 5$ nonnegative bisymmetric matrix with the eigenvalues $\lambda_1, \lambda_2, \lambda_3, \lambda_4, \lambda_5$ is the matrix

$$L_1 = \begin{bmatrix}
\dfrac{\lambda_1+\lambda_5}{2}& &&&\dfrac{\lambda_1-\lambda_5}{2}\\
&\dfrac{\lambda_2+\lambda_4}{2}&&\dfrac{\lambda_2-\lambda_4}{2}&\\
&&\lambda_3&&\\
&\dfrac{\lambda_2-\lambda_4}{2}&&\dfrac{\lambda_2+\lambda_4}{2}&\\
\dfrac{\lambda_1-\lambda_5}{2}&&&&\dfrac{\lambda_1+\lambda_5}{2}
\end{bmatrix}.$$

If $\lambda_1 \geq 0 > \lambda_2 \geq \lambda_3 \geq \lambda_4 \geq \lambda_5$ then a $5 \times 5$ nonnegative bisymmetric matrix with the eigenvalues $\lambda_1, \lambda_2, \lambda_3, \lambda_4, \lambda_5$ is the matrix

$$L_2=\begin{bmatrix}
0  &a   &b    &a             &-\lambda_5\\
a  &0   &c    &-\lambda_4    &a\\
b  &c   &\sum_{i=1}^5 \lambda_i&   c   &b\\
a  &-\lambda_4&c    &0       &a\\
-\lambda_5&a&b&a&0
\end{bmatrix},$$
where 
$$a=\dfrac{1}{2}\sqrt{\dfrac{(\lambda_1+\lambda_5)(\lambda_2+\lambda_5)(\lambda_3+\lambda_4)}{\lambda_1+\lambda_2-\lambda_3+\lambda_5}},$$ 

$$b=\sqrt{\dfrac{-(\lambda_1+\lambda_5)(\lambda_2+\lambda_5)(\lambda_1+\lambda_2+\lambda_4+\lambda_5)}{2(\lambda_1+\lambda_2-\lambda_3+\lambda_5)}},$$ 

\noindent and

$$c=\sqrt{\dfrac{-(\lambda_3+\lambda_4)(\lambda_1+\lambda_2+\lambda_4+\lambda_5)}{2}}.$$

Now we consider the case $\lambda_1 \geq \lambda_2 \geq 0 > \lambda_3 \geq \lambda_4 \geq \lambda_5$. If $\lambda_1+ \lambda_3+\lambda_4 \geq 0$ and $\lambda_2+\lambda_5 < 0$ then the matrix $L_2$ in the previous case is nonnegative and hence it is our desired matrix. If $\lambda_2+\lambda_5 \geq 0$ and $\lambda_1+\lambda_3+\lambda_4 \geq 0$ then a $5 \times 5$ nonnegative bisymmetric matrix with the eigenvalues $\lambda_1, \lambda_2, \lambda_3, \lambda_4, \lambda_5$ is the matrix

$$L_3=\begin{bmatrix}
\dfrac{\lambda_2+\lambda_5}{2}&0&0&0&\dfrac{\lambda_2-\lambda_5}{2}\\
0&0&d&-\lambda_4&0\\
0&d&\lambda_1+\lambda_3+\lambda_4&d&0\\
0&-\lambda_4&d&0&0\\
\dfrac{\lambda_2-\lambda_5}{2}&0&0&0&\dfrac{\lambda_2+\lambda_5}{2}
\end{bmatrix},$$\\
where $d=\sqrt{\dfrac{-(\lambda_3+\lambda_4)(\lambda_1+\lambda_4)}{2}}$. Finally, if $\lambda_2+\lambda_5 \geq 0$ and $\lambda_1+\lambda_3+\lambda_4 < 0$ then $\sum_{i=1}^5 \lambda_i < \lambda_2+\lambda_5$ contradicting a necessary condition of the result (Lemma $1$ in \cite{loewy}) by R. Loewy and J. J. McDonald and hence in this case there is no nonnegative bisymmetric matrix with the eigenvalues $\lambda_1, \lambda_2, \lambda_3, \lambda_4, \lambda_5$.

Next, we consider the case $\lambda_1 \geq  \lambda_2 \geq \lambda_3 \geq 0 > \lambda_4 \geq \lambda_5.$ If $\lambda_1+\lambda_2+\lambda_4+\lambda_5 \geq 0$ and $\lambda_2+\lambda_5 \geq 0$, then $\lambda_2+\lambda_4 \geq 0$ and the matrix $L_1$ is our desired matrix.
If $\lambda_1+\lambda_2+\lambda_4+\lambda_5 \geq 0$, $\lambda_2+\lambda_5 < 0$ and $\lambda_3+\lambda_4 \geq 0,$ then the matrix 
$$L_4=\begin{bmatrix}
\dfrac{\lambda_3+\lambda_4}{2}&0&0&0&\dfrac{\lambda_3-\lambda_4}{2}\\
0&0&g&-\lambda_5&0\\
0&g&\lambda_1+\lambda_2+\lambda_5&g&0\\
0&-\lambda_5&g&0&0\\
\dfrac{\lambda_3-\lambda_4}{2}&0&0&0&\dfrac{\lambda_3+\lambda_4}{2}
\end{bmatrix},$$
\noindent where $g=\sqrt{\dfrac{-(\lambda_2+\lambda_5)(\lambda_1+\lambda_5)}{2}}$, \ will be our desired matrix.
If $\lambda_1+\lambda_2+\lambda_4+\lambda_5 \geq 0$, $\lambda_2+\lambda_5 < 0$ and $\lambda_3+\lambda_4 < 0$, then again the matrix $L_2$ can be used. We do not have a complete solution for the case when $\lambda_1+\lambda_2+\lambda_4+\lambda_5 < 0$. However, some partial result (Corollary \ref{cor5by5}) will be given at the end of Section $3$.
\bigskip

Summing up what we have done so far, we get the following result.
\bigskip
\begin{prop} \label{prop5by5}
Let $\lambda_1 \geq \lambda_2 \geq \lambda_3 \geq \lambda_4 \geq \lambda_5 \geq -\lambda_1$ be real numbers such that $\sum_{i=1}^5 \lambda_i \geq 0$. If $\lambda_1, \lambda_2, \lambda_3,\lambda_4, \lambda_5$ satisfies one of the following conditions
\begin{itemize}

\item[(1)] $\lambda_4 \geq 0$ or $\lambda_2 \leq 0$,
\item[(2)] $\lambda_2 \geq 0 > \lambda_3$ and $\sum_{i=1}^5 \lambda_i \geq \lambda_2+\lambda_5$,
\item[(3)] $\lambda_3 \geq 0 > \lambda_4$ and $\lambda_1 + \lambda_2 +\lambda_4+\lambda_5 \geq 0$,
\end{itemize}
then there is a $5 \times 5$ nonnegative bisymmetric matrix with the eigenvalues $\lambda_1$, $\lambda_2$, $\lambda_3$, $\lambda_4$, $\lambda_5.$
\end{prop}
\bigskip
Proposition \ref{prop5by5} will be employed to obtain our main result, Theorem \ref{main5by5} in Section $3$.
\bigskip
\section{A necessary and sufficient condition of the BNIEP for $5 \times 5$ traceless matrices}
\bigskip
In this section, we give a characterization of bisymmetric nonnegative traceless $5 \times 5$ matrices according to the characterization of symmetric nonnegative traceless $5 \times 5$ matrices given by O. Spector in \cite{oren}. To obtain our main result, we use Proposition \ref{prop5by5} in the previous section along with a useful result of Cantoni and Butler in \cite{can} which gives a standard form of a bisymmetric matrix. This result of Cantoni and Butler (Lemma 2(ii) and Lemma 3(b) in \cite{can}) states that a $(2m+1) \times (2m+1)$ matrix $Q$ is bisymmetric if and only if $Q$ is of the form $\begin{bmatrix}
A   & x  & JCJ\\
x^T & p  & x^TJ\\
C   & Jx & JAJ
\end{bmatrix},$ 
where $A$ and $C$ are $m \times m$ matrices, $A=A^T$, $JC^TJ=C$, $x$ is an $m \times 1$ matrix and $p$ is a real number, and in this case the eigenvalue of $Q$ can be partitioned into the eigenvalue of  $A-JC$ and the eigenvalues of 
$\begin{bmatrix}
p          & \sqrt{2} x^T \\
\sqrt{2} x & A+JC
\end{bmatrix}$.

So, if $p \in \mathbb{R}$, $x$ is an $m \times 1$ matrix and $X$, $Y$ are $m \times m$ symmetric matrices, then $A=\dfrac{X+Y}{2}$ and $JC=\dfrac{X-Y}{2}$ are symmetric and hence the matrix
$$Q=\dfrac{1}{2}\begin{bmatrix}
X+Y          & \sqrt{2}x    & (X-Y)J \\
\sqrt{2}x^T  & p	        & \sqrt{2}x^TJ \\
J(X-Y)       & \sqrt{2}Jx   & J(X+Y)J
\end{bmatrix}$$ is a bisymmetric matrix with the eigenvalues obtained from the matrices $\begin{bmatrix}
p & x^T \\
x & X
\end{bmatrix}$ and $Y$.
In terms of $A$ and $C$, the matrix $Q$ is of the form 
$\begin{bmatrix}
A                     & \dfrac{x}{\sqrt{2}}  & JCJ \\
\dfrac{x^T}{\sqrt{2}} &         p            & \dfrac{x^TJ}{\sqrt{2}} \\
C                     & \dfrac{Jx}{\sqrt{2}} & JAJ 
\end{bmatrix}$. Note that, $Q$ is a nonnegative matrix if and only if $p \geq 0$ and $x, X+Y$ and $X-Y$ are nonnegative matrices. So $Q$ is nonnegative if we have the condition that 
$\begin{bmatrix}
p & x^T \\
x & A+JC
\end{bmatrix}$, 
$A$ and $C$ are nonnegative. Therefore, to construct a $(2m+1)\times (2m+1)$ bisymmetric nonnegative matrix from $m \times m$ symmetric matrices $X$ and $Y$, $p \in \mathbb{R}$ and an $m \times 1$ matrix $x$, we need to find matrices $A$ and $C$ such that $X=\dfrac{A+JC}{2}$, $Y=\dfrac{A-JC}{2}$ and $\begin{bmatrix}
p & x^T \\
x & A+JC
\end{bmatrix}$, $A$ and $C$ are nonnegative. We will use this method in the proof of Theorem \ref{5by5unknowcase1}

Before we proof our main theorem, we first need the following auxiliary lemma.
\bigskip

\begin{lem} \label{5by5lemma}
Let $ 1 \geq \lambda_2 \geq \lambda_3 > 0 > \lambda_4 \geq \lambda_5 \geq -1$ be real numbers such that $1+\sum_{i=2}^5 \lambda_i = 0$ and \ $1+\sum_{i=2}^5 \lambda_{i}^3 \geq 0.$ If \ $1 + \lambda_2 + \lambda_4 + \lambda_5 < 0$, then 
$$\lambda_2 \lambda_4 - \lambda_3 \lambda_5 - \lambda_3 - \lambda_5 \geq 0 \text{ \ and \ } \lambda_2\lambda_4 - \lambda_3 \lambda_5 - \lambda_3 - \lambda_5 \geq \dfrac{\lambda_3 \lambda_5}{\lambda_2+\lambda_4}.$$
\begin{proof}
Suppose that $1 + \lambda_2 + \lambda_4 + \lambda_5 < 0$. Since $1 \geq \lambda_2 \geq \lambda_3 > 0 > \lambda_4 \geq \lambda_5 \geq -1$ and $1 + \lambda_2 + \lambda_4 + \lambda_5 < 0$, we have $\lambda_2 + \lambda_5 < 0.$ Since $1 > \lambda_2 \geq \lambda_3$ and $-\lambda_5 > \lambda_4 \geq -1$, by the rearrangement inequality, we have $-\lambda_5 + \lambda_2\lambda_4 - \lambda_3 \geq -\lambda_2 + \lambda_3\lambda_4 - \lambda_5.$ So, 
$$-\lambda_5 + \lambda_2\lambda_4 - \lambda_3 - \lambda_3\lambda_5 \, \ \geq -\lambda_2 + \lambda_3\lambda_4 - \lambda_5-\lambda_3\lambda_5= \lambda_3 (\lambda_4 - \lambda_5) - (\lambda_2+\lambda_5) \geq 0.$$
Therefore, we have $\lambda_2 \lambda_4 - \lambda_3 \lambda_5 - \lambda_3 - \lambda_5 \geq 0.$ 

Since $\lambda_5 = -(1+\lambda_2+\lambda_3+\lambda_4)$,
$$\displaystyle 0 \leq 1+\sum_{i=2}^5 \lambda_{i}^3= 1+\lambda_2^3+\lambda_3^3+\lambda_4^3 -(1+\lambda_2+\lambda_3+\lambda_4)^3.$$ 
\noindent Simplifying the right-hand side of the previous inequality, we obtain






$0 \leq \lambda_2\lambda_5-\lambda_2^2\lambda_4+\lambda_2\lambda_3\lambda_5+\lambda_2\lambda_3+\lambda_2\lambda_5-\lambda_2\lambda_4^2+\lambda_3\lambda_4\lambda_5+\lambda_3\lambda_4+\lambda_4\lambda_5.$

\noindent Therefore, 

\qquad  $\lambda_2\lambda_5 \geq \lambda_2^2\lambda_4-\lambda_2\lambda_3\lambda_5-\lambda_2\lambda_3-\lambda_2\lambda_5+\lambda_2\lambda_4^2-\lambda_3\lambda_4\lambda_5-\lambda_3\lambda_4-\lambda_4\lambda_5$

\qquad \qquad \,\ $= (\lambda_2\lambda_4-\lambda_3\lambda_5-\lambda_3-\lambda_5)(\lambda_2+\lambda_4).$

\noindent Since \ $\lambda_2+\lambda_4<0$, we have $\dfrac{\lambda_3 \lambda_5}{\lambda_2+\lambda_4}$ \ $\leq \ \lambda_2\lambda_4-\lambda_3\lambda_5-\lambda_3-\lambda_5$.
\end{proof}
\end{lem}

Now, we are ready to prove the following result.
\medskip

\begin{thm} \label{5by5unknowcase1}
Let $\lambda_1 \geq \lambda_2 \geq \lambda_3 > 0 > \lambda_4 \geq \lambda_5 \geq -\lambda_1$ be real numbers such that $\sum_{i=1}^5 \lambda_i = 0$ and $\lambda_1+\lambda_2+\lambda_4+\lambda_5 <0.$ Then $\lambda_1, \lambda_2, \lambda_3, \lambda_4, \lambda_5$ is the eigenvalue of a $5 \times 5$ nonnegative bisymmetric matrix if and only if \ $\sum_{i=1}^5 \lambda_{i}^3 \ \geq  0.$
\begin{proof}
\noindent Normalizing by the Perron root, without loss of generality, we assume that $\lambda_1 = 1$.

If there is a $5 \times 5$ nonnegative bisymmetric matrix $Q$ with the eigenvalues $1, \lambda_2, \lambda_3, \lambda_4, \lambda_5$, then $Q^3$ is a nonnegative diagonalizable matrix and hence $\displaystyle 1+\sum_{i=2}^5 \lambda_{i}^{3} =$ Tr($Q^3$) $\geq 0$. So $\displaystyle 1+\sum_{i=2}^5 \lambda_{i}^{3} \geq 0$ is a necessary condition.

Assume that 
$\displaystyle 1+\sum_{i=2}^5 \lambda_{i}^3 \geq 0.$ Then, by Lemma \ref{5by5lemma}, we have $\lambda_2\lambda_4 - \lambda_3\lambda_5 - \lambda_3 - \lambda_5 \geq 0$ and $\lambda_2\lambda_4 - \lambda_3 \lambda_5 - \lambda_3 - \lambda_5 \geq \dfrac{\lambda_3 \lambda_5}{\lambda_2+\lambda_4}.$ Let $A$ and $C$ be nonnegative matrices such that
$$A-JC=\begin{bmatrix}
\lambda_2+\lambda_4 & \sqrt{-\lambda_2\lambda_4}\\
\sqrt{-\lambda_2\lambda_4}&0
\end{bmatrix} 
\text{ \ and \ }
A+JC =\begin{bmatrix}
-(\lambda_2+\lambda_4) & \sqrt{-\lambda_2\lambda_4}\\
\sqrt{-\lambda_2\lambda_4}&0
\end{bmatrix}.$$
Then $A=\begin{bmatrix}
0                          &  \sqrt{-\lambda_2\lambda_4} \\
\sqrt{-\lambda_2\lambda_4} & 0
\end{bmatrix}$ and 
$C=\begin{bmatrix}
           0               &  0 \\
-(\lambda_2+\lambda_4)     & 0
\end{bmatrix}$. Note that $A-JC$ has \ $\lambda_2$ and $\lambda_4$ \ as its eigenvalues. Now, let
$$U_{a,b}=\begin{bmatrix}
0                      &                 a           & b \\
a                      & -(\lambda_2+\lambda_4)      & \sqrt{-\lambda_2\lambda_4} \\
b                      &  \sqrt{-\lambda_2\lambda_4} & 0 \\
\end{bmatrix},$$ 

\noindent where $a, b$ are nonnegative real numbers. We want to show that there are real numbers $a$ and $b$ such that the matrix $U_{a,b}$ has the eigenvalues $1, \lambda_3, \lambda_5$.

Suppose that $1, \lambda_3, \lambda_5$ are the eigenvalues of $U_{a,b}.$ Then by comparing the coefficients of $(x-1)(x-\lambda_3)(x-\lambda_5)$ and the characteristic polynomial of $U_{a,b}$ we have the system of equations
$$a^2+b^2=\lambda_2\lambda_4-\lambda_3\lambda_5-\lambda_3-\lambda_5 \text{ \ and \ } 
-(\lambda_2+\lambda_4)b^2 - 2ab\sqrt{-\lambda_2\lambda_4}=-\lambda_3\lambda_5.$$

\noindent We observe that the first equation represent a circle with radius 
$$r = \sqrt{\lambda_2\lambda_4-\lambda_3\lambda_5-\lambda_3-\lambda_5}$$ and the second equation gives the hyperbola with $b$-intercept at $b= \pm  \sqrt{\dfrac{\lambda_3\lambda_5}{\lambda_2+\lambda_4}}$ in $ab$-coordinates. By Lemma \ref{5by5lemma}, $r^2 \geq \dfrac{\lambda_3 \lambda_5}{\lambda_2+\lambda_4}$. So, the hyperbola intersects $b$-axis inside the circle. Using the quadratic formula to the equation
$-(\lambda_2+\lambda_4)b^2 - 2ab\sqrt{-\lambda_2\lambda_4}=-\lambda_3\lambda_5$ we have

 $$b=\dfrac{2a\sqrt{-\lambda_2\lambda_4}\pm\sqrt{-4a^2\lambda_2\lambda_4+4(\lambda_2+\lambda_4)\lambda_3\lambda_5}}{-2(\lambda_2+\lambda_4)}.$$
 
We see that the function $f$ defined by $$f(a)=\dfrac{2a\sqrt{-\lambda_2\lambda_4}+\sqrt{-4a^2\lambda_2\lambda_4+4(\lambda_2+\lambda_4)\lambda_3\lambda_5}}{-2(\lambda_2+\lambda_4)}$$ 

\noindent for $a \geq 0$, is increasing and $f(a) \geq 0$ for $a \geq 0$. Therefore the hyperbola and the circle are intersect at some point in the first quadrant. So, there are nonnegative real numbers $a_0$ and $b_0$ satisfying our system and hence the matrix $U_{a_0,b_0}$ has the eigenvalue $1, \lambda_3, \lambda_5$. Finally, we define the matrix 
$Q=\begin{bmatrix}
A&\dfrac{x}{\sqrt{2}}&JCJ\\
\dfrac{x^T}{\sqrt{2}}&0&\dfrac{x^TJ}{\sqrt{2}}\\
C&\dfrac{Jx}{\sqrt{2}}&JAJ
\end{bmatrix}$ where $x^T = \begin{bmatrix}
a_0 & b_0
\end{bmatrix}.$ Then $Q$ is our desired matrix.
\end{proof}
\end{thm}
\bigskip
Now, we state and prove our main result for this section.
\bigskip
\begin{thm} \label{main5by5}
Let $\lambda_1 \geq  \lambda_2 \geq \lambda_3 \geq \lambda_4 \geq \lambda_5 \geq -\lambda_1$ be real numbers such that $\sum_{i=1}^5 \lambda_i = 0$. Then $\lambda_1, \lambda_2,\lambda_3,\lambda_4,\lambda_5$ is the eigenvalue of a $5 \times 5$ nonnegative bisymmetric matrix if and only if 
$$\lambda_2 +\lambda_5 \leq 0 \text{ \ and \ } \sum_{i=1}^5 \lambda_i^3 \geq 0.$$

\begin{proof}
We have seen that $\sum_{i=1}^5 \lambda_i^3 \geq 0$ is a necessary condition. In \cite{loewy}, R. Loewy and J. J. McDonald proved that for any $5 \times 5$ nonnegative symmetric matrix $A$  with the eigenvalues $\alpha_1 \geq \alpha_2 \geq \alpha_3 \geq \alpha_4 \geq \alpha_5$, Tr($A$) $\geq \alpha_2+\alpha_5$. So, for any traceless bisymmetric matrix $Q$ with the eigenvalues $\alpha_1 \geq \alpha_2 \geq \alpha_3 \geq \alpha_4 \geq \alpha_5$, we have $\alpha_2+\alpha_5 \leq$ Tr($Q$) $=0$. So, we prove the direct part of our assertion.

Conversely, assume that $\lambda_2+\lambda_5 < 0$ and $\sum_{i=1}^5 \lambda_i^3 \geq 0$.

Case 1: If $\lambda_5 \geq 0$, then all $\lambda_i$'s are $0$ and the zero $5 \times 5$ matrix is our desired matrix.

Case 2: If $\lambda_4 \geq 0 > \lambda_5$ or $\lambda_2 \geq 0 > \lambda_3$ or $\lambda_1 \geq 0 > \lambda_2$, then we can apply Propositon \ref{prop5by5}.

Case 3: Suppose that $\lambda_3 \geq 0 > \lambda_4$. If $\lambda_1+\lambda_2+\lambda_4+\lambda_5 \geq 0,$ then we can use Proposition \ref{prop5by5}. If $\lambda_1+\lambda_2+\lambda_4+\lambda_5 < 0,$ then $\lambda_3 \neq 0$ and hence Theorem \ref{5by5unknowcase1} can be applied.
\end{proof}
\end{thm}

Finally, we give a partial result concerning the case that we do not have a complete answer, that is the case of finding a $5\times5$ bisymmetric nonnegative matrix having the eigenvalues $\lambda_1, \lambda_2, \lambda_3, \lambda_4, \lambda_5$ satisfying the condition $\lambda_1 \geq  \lambda_2 \geq \lambda_3 \geq 0 > \lambda_4 \geq \lambda_5 \geq -\lambda_1$, $\sum_{i=1}^5 \lambda_i > 0$ and $\lambda_1+\lambda_2+\lambda_4+\lambda_5 <0.$ Note that, in this condition, we have $\lambda_1 \neq 0$ and $\lambda_3 \neq 0.$ Therefore, without loss of generality, we may assume that $\lambda_1=1$ in this condition.
\bigskip
\begin{cor} \label{cor5by5}
Let $1 \geq  \lambda_2 \geq \lambda_3 > 0 > \lambda_4 \geq \lambda_5 \geq -1$ be real numbers such that $1+\sum_{i=2}^5 \lambda_i > 0$ and \ $1+\lambda_2+\lambda_4+\lambda_5 < 0$. If \ 
$\lambda_2\lambda_4 - \lambda_3 \lambda_5 - \lambda_3 - \lambda_5 \geq \dfrac{-\lambda_3 \lambda_5}{1+\lambda_3+\lambda_5},$ then there is a $5 \times 5$ nonnegative bisymmetric matrix with the eigenvalues $1,  \lambda_2, \lambda_3, \lambda_4, \lambda_5.$
\begin{proof}
Again, by setting  
$$A-JC=\begin{bmatrix}
\lambda_2+\lambda_4 & \sqrt{-\lambda_2\lambda_4}\\
\sqrt{-\lambda_2\lambda_4}&0
\end{bmatrix}, \ 
A+JC =\begin{bmatrix}
1+\lambda_3+\lambda_5 & \sqrt{-\lambda_2\lambda_4}\\
\sqrt{-\lambda_2\lambda_4}&0
\end{bmatrix}$$ 

\noindent and $U_{a,b}=\begin{bmatrix}
0                     &                 a           &     b        \\
a                     &    1+\lambda_3+\lambda_5    & \sqrt{-\lambda_2\lambda_4} \\
b                     &  \sqrt{-\lambda_2\lambda_4} &    0          \\
\end{bmatrix}$ and then comparing the coefficient of characteristic polynomials of $U_{a,b}$ and the polynomials $(x-1)(x-\lambda_3)(x-\lambda_5)$, we have the system
$$a^2+b^2=\lambda_2\lambda_4-\lambda_3\lambda_5-\lambda_3-\lambda_5 \text{ \ and \ }
(1+\lambda_3+\lambda_5)b^2 - 2ab\sqrt{-\lambda_2\lambda_4}=-\lambda_3\lambda_5.$$
Using the idea in the proof of Theorem \ref{5by5unknowcase1}, the inequality condition $\lambda_2\lambda_4 - \lambda_3 \lambda_5 - \lambda_3 - \lambda_5 \geq \dfrac{-\lambda_3 \lambda_5}{1+\lambda_3+\lambda_5}$ give us a nonnegative solution $(a_0, b_0)$ of the above system. Then the matrix 

$$Q=\begin{bmatrix}
\dfrac{1+\lambda_2+\lambda_3+\lambda_4+\lambda_5}{2} & \sqrt{-\lambda_2\lambda_4} & \dfrac{a_0}{\sqrt{2}} & 0 & \dfrac{1+\lambda_3+\lambda_5-\lambda_2-\lambda_4}{2} \\
\sqrt{-\lambda_2\lambda_4} & 0 & \dfrac{b_0}{\sqrt{2}} & 0 & 0\\
\dfrac{a_0}{\sqrt{2}}  &  \dfrac{b_0}{\sqrt{2}} & 0 & \dfrac{b_0}{\sqrt{2}} & \dfrac{a_0}{\sqrt{2}} \\
0   & 0 &\dfrac{b_0}{\sqrt{2}}  & 0 & \sqrt{-\lambda_2\lambda_4}\\
\dfrac{1+\lambda_3+\lambda_5-\lambda_2-\lambda_4}{2} & 0 &\dfrac{a_0}{\sqrt{2}}  & \sqrt{-\lambda_2\lambda_4} & \dfrac{1+\lambda_2+\lambda_3+\lambda_4+\lambda_5}{2}
\end{bmatrix}$$
is our desired matrix.
\end{proof}
\end{cor}
\medskip
\begin{exa}
The set of real numbers $\lbrace 1, 0.3, 0.2, -0.7, -0.8 \rbrace$ satisfies the condition of Theorem \ref{5by5unknowcase1}. Using the construction in the proof of Theorem \ref{5by5unknowcase1}, we let $A$ and $C$ be two matrices such that 
$$A-JC=\begin{bmatrix}
-0.4        & \sqrt{0.21}\\
\sqrt{0.21} & 0
\end{bmatrix} \text{ \ and \ } A+JC=\begin{bmatrix}
0.4		    & \sqrt{0.21} \\
\sqrt{0.21} & 0
\end{bmatrix} $$
and let $U_{a,b}=\begin{bmatrix}
0   &    a            &     b \\
a   &    0.4          &  \sqrt{0.21} \\
b   &    \sqrt{0.21}  &     0
\end{bmatrix}$ be such that $U_{a, b}$ has the eigenvalues $1$, $0.2$, $-0.8$. That is, we want to find nonnegative real numbers $a$ and $b$ satisfying the system
$$a^2+b^2 = 0.55 \text{ \ and \ } 0.4b^2-2\sqrt{0.21}ab = 0.16.$$ 
One nonnegative solution $(a_0, b_0)$ of this system is given by $a_0=\sqrt{\dfrac{51-3\sqrt{273}}{200}}$ and $b_0=\sqrt{\dfrac{59+3\sqrt{273}}{200}}$. Therefore the matrix $$Q=\begin{bmatrix}
A              & \dfrac{x}{2}  & JCJ \\
\dfrac{x^T}{2} & 0             & \dfrac{x^T}{2}J \\
C              & J\dfrac{x}{2} & JAJ
\end{bmatrix}=
\begin{bmatrix}
0           & \sqrt{0.21} & \dfrac{a_0}{\sqrt{2}} & 0           &    0.4 \\
\sqrt{0.21} & 0           & \dfrac{b_0}{\sqrt{2}} & 0           &      0   \\
\dfrac{a_0}{\sqrt{2}}     & \dfrac{b_0}{\sqrt{2}} & 0           & \dfrac{b_0}{\sqrt{2}} & \dfrac{a_0}{\sqrt{2}}\\
0           & 0           & \dfrac{b_0}{\sqrt{2}} & 0           & \sqrt{0.21} \\
0.4         & 0           & \dfrac{a_0}{\sqrt{2}} & \sqrt{0.21} & 0
\end{bmatrix},$$ 
where $x^T=\begin{bmatrix}
a_0&b_0
\end{bmatrix}$, is a nonnegative bisymmetric matrix with the eigenvalues $1, 0.3, 0.2, -0.7, -0.8$.
\end{exa}

\end{document}